# The parametric continuation method for determining steady state diagrams

Marek Berezowski

Silesian University of Technology, Faculty of Applied Mathematics, Gliwice, Poland

E-mail: marek.berezowski@polsl.pl

## Abstract

The paper discusses a simple method of using the parametric continuation method to designate complex diagrams of steady states. The main advantage of the discussed approach is the fact that it does not require the installation of huge professional IT systems. The deliberations are illustrated by examples of a reactor models described by the algebraic and differential equations. Both models render solutions in the form of the so called multiple steady states.

## 1. Introduction

Mathematical models that describe equipment offer ambiguous solutions, due to the so called: "multiple steady states phenomenon", where specifically determined parameter values correspond to more than one variable values (Berezowski, 1987, Berezowski and Burghardt, 1989, Burghardt and Berezowski, 1990 and Subramanian and Balakotaiah, 1996). The steady state diagrams created on the bases of this may have a very complicated nature, presenting a tangle of lines on the plane or in space (Berezowski, 2000). In view of a non-linear character of the models, often described by means of differential equations, their analytical solution is impossible, thus, leaving space for numerical methods. Nevertheless, the analytical methods also entail essential calculation problems arising from the ambiguity of solutions and small distances between specific steady states (Berezowski, 2000—Fig. 3). In such cases, iteration methods that render solutions dependent on the initial values are practically useless. Therefore, the parametric continuation method, which is not based on numerical iteration, may be helpful, especially that it is not sensitive to the degree of complexity of the graph (Seydel, 1994 and Shalashin and Kuznetsov, 2003).

In this paper simple ways of the application of the parametric continuation method are discussed on the bases of reactor models. The first example involves a tank reactor described by simple differential equations, whereas the second one concerns tubular reactor, described by simple differential equations of the second order with boundary conditions. Unlike the manner based on the so called: " fictitious dynamics" (Berezowski, 2000) the discussed method does not require the installation of huge professional IT systems (Doedel et al., 1997).

## 2. General description of determining the steady states by means of the parametric continuation method

Let us assume an *n*-dimensional model:

$$\overline{F}(\overline{y},p) = \overline{0} \tag{1}$$

from which *n* characteristics (diagrams) should be designated as $y_i(p)$ in a specified interval of the variability of parameter *p*. If the non-linear equation system (1) is mathematically confounding, the analytical designation of particular diagrams seems to be impossible. In such a case, it is safe to fall back on numerical methods, one of which could involve an iteration solution of the system for the successive values of *p*, for example, on the bases of Newton's method. Yet, the main problem arises when the solutions of system (1) are ambiguous, i.e. when a given value of *p* corresponds to more than one vector $\overline{y}(p)$. In such a case, the sought solutions are very sensitive to changes in the initial values, especially when the solutions are located in very close proximity. In the majority of cases, it is practically impossible to derive correct diagrams. Thus, other methods have to be employed, for example, the parametric continuation method, the description and application of which is presented below.

By designating the total differential of system (1) we obtain

$$\overline{\overline{J}}\, d\overline{y} + \overline{w}\, dp = \overline{0} \tag{2}$$

where $\overline{\overline{J}}$ is Jacobi matrix in the following form:

$$\overline{\overline{J}} = \left\{\frac{\partial f_i}{\partial y_j}\right\};\quad i = 1, n \tag{3}$$

whereas, $\overline{w}$ is the partial derivative vector calculated in terms of parameter *p*:

$$\overline{w} = \left\{\frac{\partial f_i}{\partial p}\right\};\quad i = 1, n \tag{4}$$

From Eq. (2) we derive that

$$d\overline{y} = -\overline{\overline{J}}^{-1}\, \overline{w}\, dp \tag{5}$$

If the above dependence is to be used in numerical calculations, differential increases $d\overline{y}$ and *dp* should be substituted by the following difference increases:

$$\Delta\overline{y} = \overline{y}_{K+1} - \overline{y}_k,\quad \Delta p = p_{k+1} - p_k \tag{6}$$

which leads to the following relation:

$$\overline{y}_{k+1} = \overline{y}_k - \overline{\overline{J}}_k^{-1}\, \overline{w}_k\, \Delta p. \tag{7}$$

If the computational process is started with precisely determined values of $\bar{y}_0$ and $p_0$ (that is, fulfilling equation (1) with big accuracy), the successive values of vector $\bar{y}_k$, determining the next points of the diagrams, are designated without the necessity of using any other iterations. However, the discussed method does not eliminate the problem involved in the ambiguity of the solutions of equations system (1). As shown in conceptual Scheme 1, at bifurcation limit points *LP* there is a change of the direction of designating parameter *p*. Thus, the points should be determined numerically in order to change the sign of increase, $\Delta p$. It should be pointed out that at points *LP*, due to the occurrence of the extreme value of parameter *p* in relation to variable *y*, differential increase $dp=0$ and matrix equation (2) is transformed as follows:

$$\bar{\bar{J}}\,d\bar{y} = \bar{0}. \qquad (8)$$

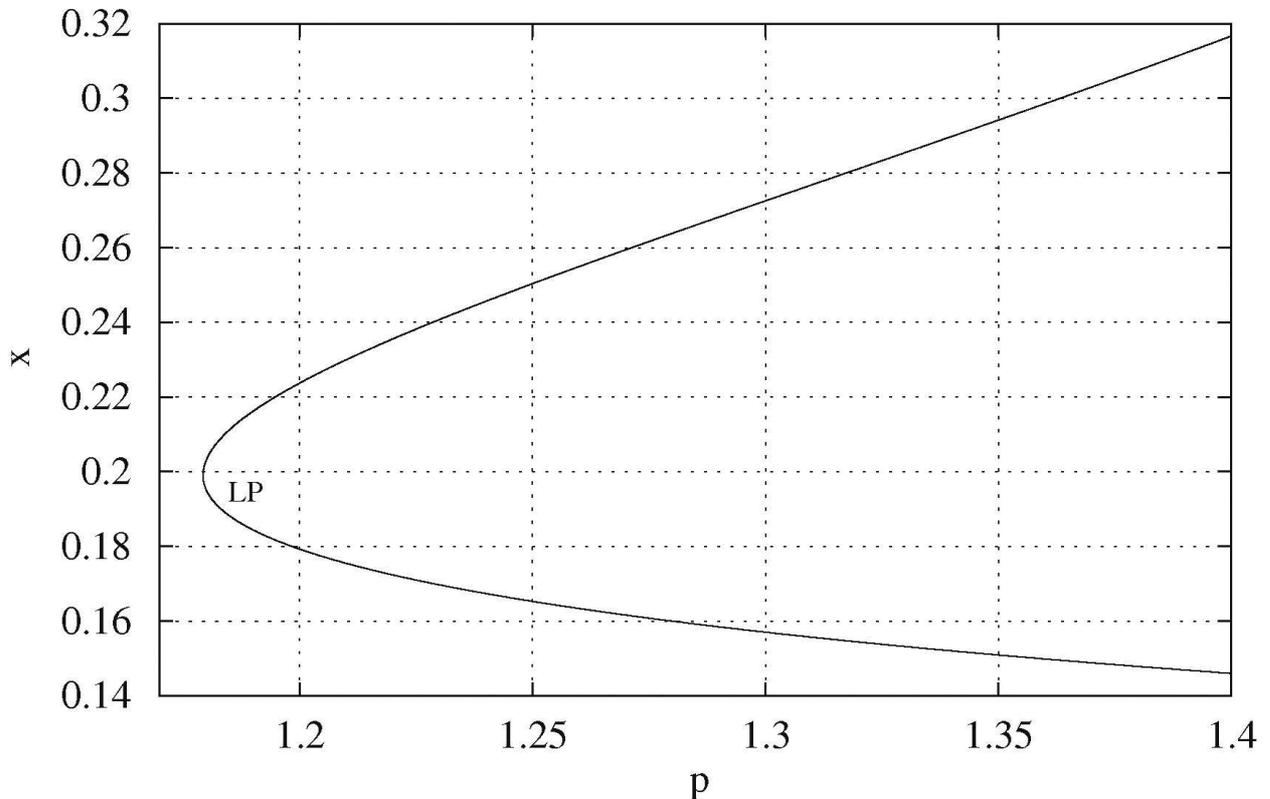

Scheme 1. Conceptual diagram.

The above equation is met only when matrix determinant $\bar{\bar{J}}$ is equal to zero:

$$\det \bar{\bar{J}} = 0 \qquad (9)$$

This means that, the determinant changes its sign at points *LP*. Hence, specific values of parameter *p* should be derived from the following relation:

$$p_{k+1} = p_k + sign(\det \bar{\bar{J}}_k)\,\Delta p \qquad (10)$$

Also in equation (7) $\Delta p$ must be multiplied by $sign(\det \bar{\bar{J}}_k)$.

At the same time, it should be emphasized that in view of the fact that the determinant zeros points $LP$, matrix $\overline{\overline{J}}$ becomes peculiar, which means that its inverse form cannot be derived. However, in view of this, there is a limit:

$$\left| \lim_{\Delta p \to 0} \frac{\Delta p}{\det \overline{\overline{J}}} \right|_{LP} < \infty \tag{11}$$

The computational process may be continued at points $LP$.

# 3. Examples of designating steady-state diagrams by means of the parametric continuation method

## 3.1. Dimensionless model of a tank reactor

The mass balance:

$$\frac{d\alpha}{d\tau} + \alpha = \phi(\alpha,\Theta) \tag{12}$$

the heat balance:

$$Le\frac{d\Theta}{d\tau} + \Theta = \phi(\alpha,\Theta) + \delta(\Theta_C - \Theta) \tag{13}$$

the kinetics function:

$$\phi(\alpha,\Theta) = Da(1-\alpha)^n \exp\left(\gamma\frac{\beta\Theta}{1+\beta\Theta}\right) \tag{14}$$

According to the definition of the steady states, they may be designated by equating both time derivatives in the above equations to zero. Accordingly, we obtain the following system of equations:

$$f_1 = -\alpha + \phi(\alpha,\Theta) = 0 \tag{15}$$

$$f_2 = -\Theta + \phi(\alpha,\Theta) + \delta(\Theta_C - \Theta) = 0 \tag{16}$$

from which $\alpha(p)$ and $\Theta(p)$ steady state characteristics may be derived for any model parameter. In such a case, the elements of Jacobi matrix assume the following form:

$$j_{11} = -1 - nDa(1-\alpha)^{n-1} \exp\left(\gamma\frac{\beta\Theta}{1+\beta\Theta}\right) \tag{17}$$

$$j_{12} = Da(1-\alpha)^n \exp\left(\gamma\frac{\beta\Theta}{1+\beta\Theta}\right)\frac{\gamma\beta}{(1+\beta\Theta)^2} \tag{18}$$

$$j_{21} = -nDa(1-\alpha)^{n-1} \exp\left(\gamma\frac{\beta\Theta}{1+\beta\Theta}\right) \tag{19}$$

$$j_{22} = -1 + Da(1-\alpha)^n \exp\left(\gamma\frac{\beta\Theta}{1+\beta\Theta}\right)\frac{\gamma\beta}{(1+\beta\Theta)^2} - \delta \tag{20}$$

If, for example, we assume that the selected bifurcation parameter is *n*, the elements of vector $\bar{w}$ will have the following form:

$$w_1 = w_2 = Da(1-\alpha)^n \exp\left(\gamma \frac{\beta\Theta}{1+\beta\Theta}\right) \ln(1-\alpha) \qquad (21)$$

At the same time, it is worth mentioning that if elements $j_{21}$ and $j_{22}$ of matrix $\bar{\bar{J}}$ are divided by *Le*, the eigenvalues of the matrix transformed in such a manner indicate the nature of the dynamical solutions of the model. Thus, if the eigenvalues of the matrix are imaginary, we are dealing with Hopf bifurcation *HB*, determining oscillations ( Jacobsen and Berezowski, 1998). So, after considering Eqs. (17), (18), (19), (20) and (21) in (7) and (10), a diagram of the tested model was derived—see Fig. 1.

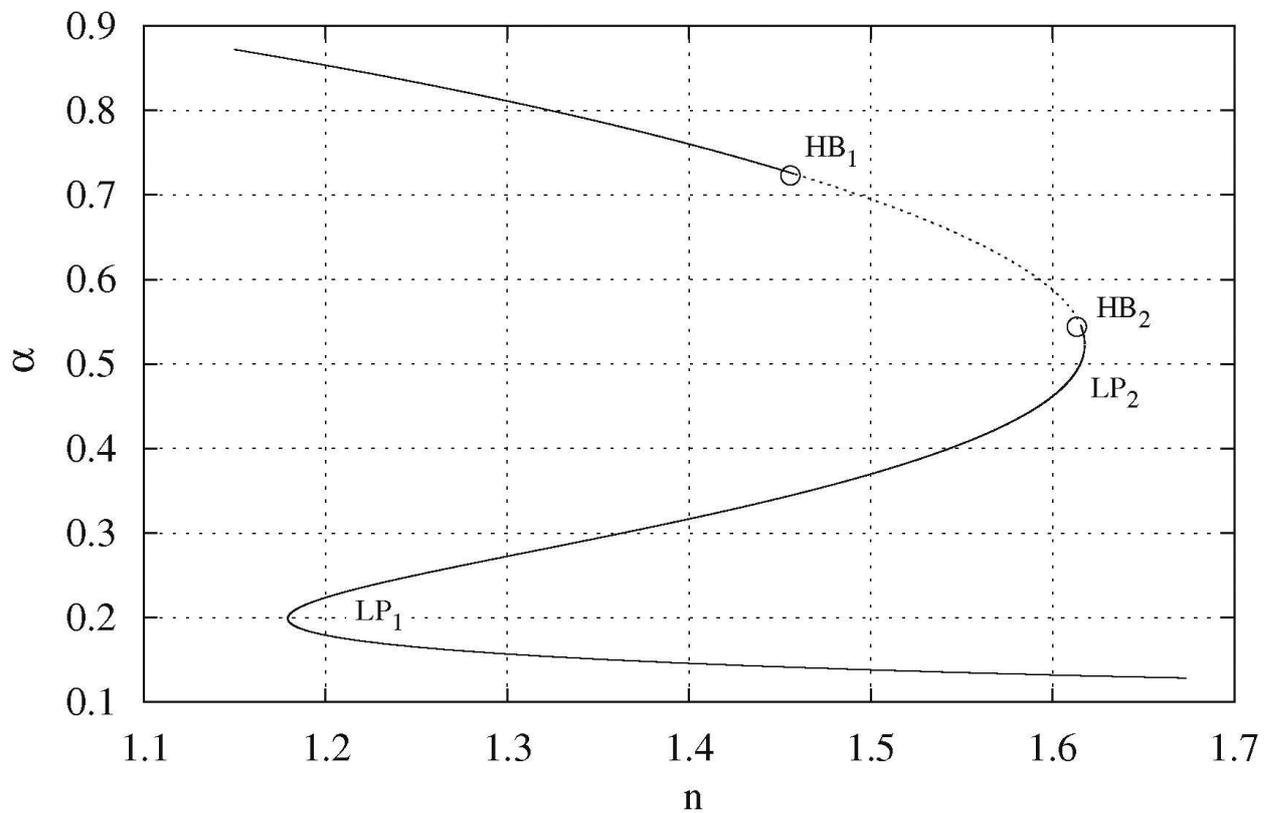

Fig. 1. Bifurcation diagram of a tank reactor.

The following parameter values were assumed in the calculations: *Le*=1.5, *Da*=0.2, *γ*=20, *β*=1, *δ*=2, *Θ*$_C$=−0.08. Bifurcation points *LP* visible in the diagram determine changes in the number of the steady states, whereas points *HB* determine the generation of oscillations.

## 3.2. Dimensionless model of an tubular reactor

The balance equations concerning the steady state have the following form:

the mass balance:

$$\frac{d\alpha}{d\xi} = \frac{1}{Pe_M}\frac{d^2\alpha}{d\xi^2} + \psi(\alpha,\Theta) \qquad (22)$$

the heat balance:

$$\frac{d\Theta}{d\xi} = \frac{1}{Pe_H}\frac{d^2\Theta}{d\xi^2} + \psi(\alpha,\Theta) \qquad (23)$$

The boundary conditions ascribed to the above system of equations are as follows:

$$\alpha(0) = \frac{1}{Pe_M}\frac{d\alpha}{d\xi}\bigg|_{\xi=0}; \quad \Theta(0) = \frac{1}{Pe_H}\frac{d\Theta}{d\xi}\bigg|_{\xi=0} \qquad (24)$$

$$\frac{d\alpha}{d\xi}\bigg|_{\xi=1} = 0; \quad \frac{d\Theta}{d\xi}\bigg|_{\xi=1} = 0 \qquad (25)$$

Assuming that $Pe_M = Pe_H = Pe$, it is easy to demonstrate that relation $\Theta = \alpha$ holds in the steady state. Under such circumstances the following system of equations is reduced to the singular equation:

$$\frac{d\alpha}{d\xi} = \frac{1}{Pe}\frac{d^2\alpha}{d\xi^2} + \phi(\alpha) \qquad (26)$$

with the following boundary conditions:

$$\alpha(0) = \frac{1}{Pe}\frac{d\alpha}{d\xi}\bigg|_{\xi=0}; \quad \frac{d\alpha}{d\xi}\bigg|_{\xi=1} = 0 \qquad (27)$$

where the kinetics function is

$$\phi(\alpha) = Da(1-\alpha)^n \exp\left(\gamma\frac{\beta\alpha}{1+\beta\alpha}\right) \qquad (28)$$

By inserting auxiliary variable:

$$u = \frac{d\alpha}{d\xi} \qquad (29)$$

the system of equations (26) and (27) is transformed to the form:

$$\frac{du}{d\xi} = Pe(u - \phi(\alpha)) \qquad (30)$$

with boundary conditions:

$$\alpha(0) = \frac{1}{Pe} u(0); \quad \alpha(1) = 0 \tag{31}$$

In the course of the analysis of the above system of equations it may be proved that the integration of equations (29) and (30) within the range from $\xi=0$ to $\xi=1$ is unstable. Thus, their numerical solution is practically impossible. Accordingly, inverse integration, i.e. within the range from $\xi=1$ to $\xi=0$ should be employed. To designate a complete diagram of the steady states the following function should be defined:

$$f(\alpha(1)) = \alpha(0) - \frac{1}{Pe} u(0) = 0 \tag{32}$$

which means, that the discussed problem is reduced to one-dimensional form. So, under such circumstances Eqs. (7) and (10) are transformed into the following form:

$$\alpha_{k+1}(1) = \alpha_k(1) - \frac{\frac{\partial f}{\partial p}\bigg|_k}{\frac{\partial f}{\partial \alpha(1)}\bigg|_k} \cdot sign\left(\frac{\partial f}{\partial \alpha(1)}\bigg|_k\right) \Delta p \tag{33}$$

$$p_{k+1} = p_k + sign\left(\frac{\partial f}{\partial \alpha(1)}\bigg|_k\right) \Delta p \tag{34}$$

The partial derivatives in the above equations should be designated numerically, which is not difficult. Assuming $Da$ as the bifurcation parameter, the diagram of the steady states shown in Fig. 2 was derived. Other assumed parameter values were $\gamma=15$, $\beta=2$, $n=1, 5$, $Pe=100$.

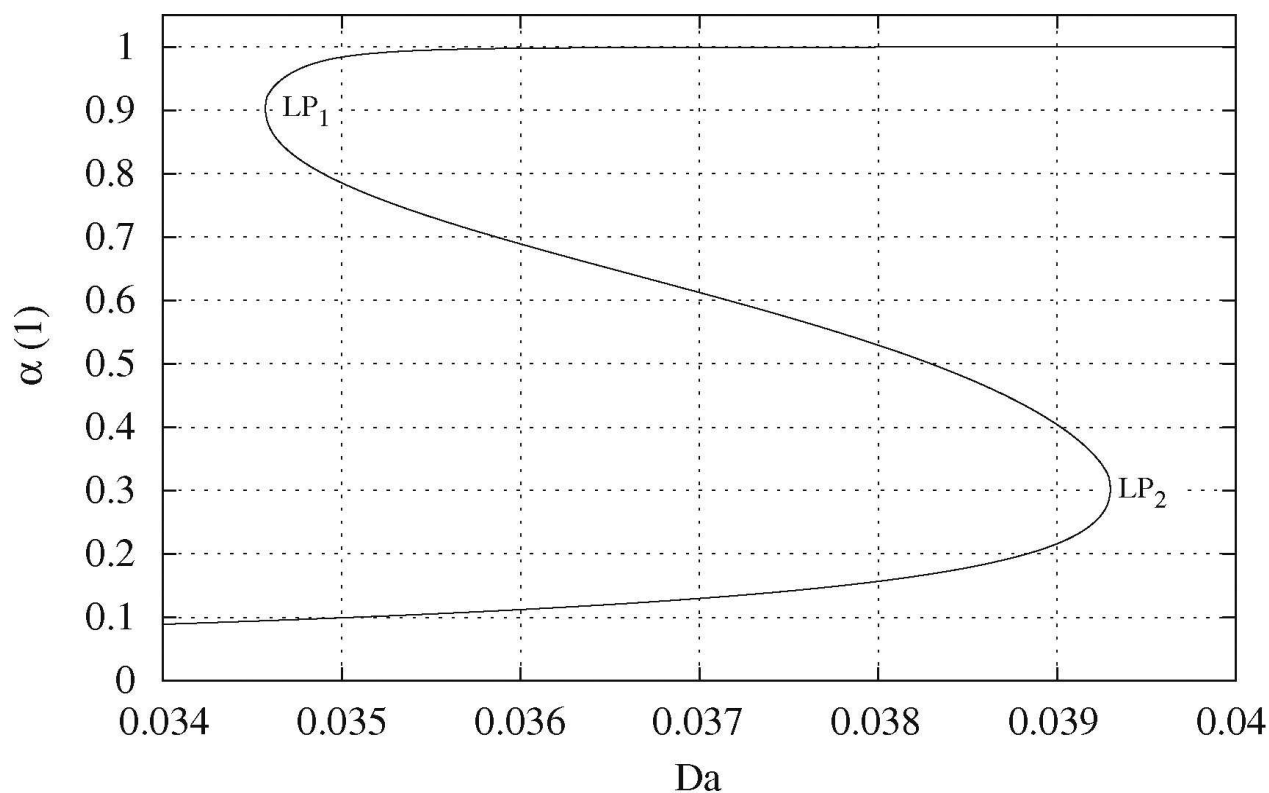
Fig. 2. Bifurcation diagram of a tubular reactor.